\documentclass[12pt,reqno]{amsart}

\usepackage{amssymb}

\usepackage{eucal}

\usepackage[all]{xy}

\font\sss=cmss8

\def\cE{{\mathcal E}}
\def\cF{{\mathcal F}}

\def\BZ{{\mathbb Z}}

\def\sD{\mbox{\sf D}}

\def\sT{\mbox{\sf T}}

\def\ssT{\mbox{\sss T}}

\def\D{\sD}

\def\Dsmall{\mbox{\sss D}}

\def\EssIm{\operatorname{Ess.\!Im}}

\def\H{\operatorname{H}}

\def\Hom{\operatorname{Hom}}

\def\Ker{\operatorname{Ker}}

\def\LTensor{\stackrel{\operatorname{L}}{\otimes}}
\def\opp{\operatorname{op}}

\def\RHom{\operatorname{RHom}}

\numberwithin{equation}{part}

\newtheorem{Lemma}{Lemma}[section]
\newtheorem{Theorem}[Lemma]{Theorem}
\newtheorem{Proposition}[Lemma]{Proposition}

\theoremstyle{definition}
\newtheorem{Definition}[Lemma]{Definition}
\newtheorem{Setup}[Lemma]{Setup}

\newtheorem{Remark}[Lemma]{Remark}

\newtheorem{Example}[Lemma]{Example}

\begin{document}

\title[Recollement for DGAs]
{Recollement for Differential Graded Algebras}

\author{Peter J\o rgensen}
\address{Department of Pure Mathematics, University of Leeds,
Leeds LS2 9JT, United Kingdom}
\email{popjoerg@maths.leeds.ac.uk}
\urladdr{http://www.maths.leeds.ac.uk/\~{ }popjoerg}




\begin{abstract} 

A recollement of triangulated ca\-te\-go\-ri\-es describes one such
category as being ``glued together'' from two others.

This paper gives a precise criterion for the existence of a
re\-col\-le\-ment of the derived category of a Differential Graded
Algebra in terms of two other such categories.

\end{abstract}

\maketitle

\setcounter{section}{-1}
\section{Introduction}
\label{sec:introduction}

\noindent
A recollement of triangulated categories is a diagram
\[
\xymatrix
{
&&\sT^{\prime} 
    \ar[rr]^{i_*} & & 
  \sT 
    \ar[rr]^{j^*}
    \ar@/_1.5pc/[ll]_{i^*} \ar@/^1.5pc/[ll]^{i^!} & &
  \sT^{\prime\prime}
    \ar@/_1.5pc/[ll]_{j_!} \ar@/^1.5pc/[ll]^{j_*} & &
}
\]
of triangulated categories and functors satisfying various conditions,
most importantly that $(i^*,i_*)$, $(i_*,i^!)$, $(j_!,j^*)$, and
$(j^*,j_*)$ are adjoint pairs (see definition \ref{dfn:recollement}
for precise details).

This notion was introduced in \cite{BBD} with the idea that $\sT$ can
be viewed as being ``glued together'' from $\sT^{\prime}$ and
$\sT^{\prime\prime}$.  The canonical example of a recollement has
$\sT$, $\sT^{\prime}$, and $\sT^{\prime\prime}$ equal to suitable
derived categories of sheaves on spaces $X$, $Z$, and $U$, where $X$
is the union of the closed subspace $Z$ and its open complement $U$.

In a more algebraic vein, an important class of triangulated
categories are the derived categories of Differential Graded modules
over Differential Graded Algebras (abbreviated below to DG modules
over DGAs).  If $R$, $S$, and $T$ are DGAs with derived categories of
left DG modules $\D(R)$, $\D(S)$, and $\D(T)$, it is therefore natural
to ask: When is there a recollement
\[
\hspace{-43ex}
\lefteqn{
\xymatrix
{
\D(S)
    \ar[rr]^{i_*} & & 
  \D(R)
    \ar[rr]^{j^*}
    \ar@/_1.5pc/[ll]_{i^*} \ar@/^1.5pc/[ll]^{i^!} & &
  \D(T) \;\; \mbox{?}
    \ar@/_1.5pc/[ll]_{j_!} \ar@/^1.5pc/[ll]^{j_*} & &
}
}
\]
The main result of this paper, theorem \ref{thm:main}, provides a
precise criterion.  In informal terms, the criterion says that
\[
  B = i_*(S) \; \mbox{ and } \; C = j_!(T)
\]
must be suitably finite objects which together generate $\D(R)$ in a
mi\-ni\-mal way.

This sheds light on several earlier results:

First, theorem \ref{thm:main} can be viewed as a generalization of the
main result of K\"{o}nig's \cite{Konig} which dealt with the situation
where $R$, $S$, and $T$ are rings.

Secondly, when $B$ and $C$ are given then $S$ and $T$ will be
constructed as endomorphism DGAs of K-projective resolutions of $B$
and $C$.  The construction of $T$ as an endomorphism DGA was
originally considered in the Morita theory developed by Dwyer and
Greenlees in \cite{DwyerGreenlees}.  I shall draw on their work, and
improve one of their results in proposition
\ref{pro:DwyerGreenleesimproved}. 

Thirdly, the construction of $T$ as an endomorphism DGA means that in
the special case $B = 0$, theorem \ref{thm:main} reduces to Keller's
theorem for DGAs, see \cite[thm.\ 4.3]{KellerDG}, which states that if
$C$ is a compact K-projective generator for $\D(R)$, then $\D(R)$ is
equivalent to $\D(\cE^{\opp})$ where $\cE$ is the endomorphism DGA of
$C$.  In the general case $B \not= 0$, theorem \ref{thm:main} can
therefore be viewed as a two object generalization of Keller's
theorem.

Finally, there is another connection to work by Keller who already in
\cite[rmk.\ 3.2]{KellerCyclic} made some remarks on recollements of
derived categories of Differential Graded Categories.

The paper is organized as follows: Section \ref{sec:embedding} gives
an embedding result for derived categories of DGAs.  Section
\ref{sec:DwyerGreenlees} recalls the Morita theory of Dwyer and
Greenlees.  And section \ref{sec:recollement} combines and develops
these themes to prove the main result.

\section{An embedding result}
\label{sec:embedding}

\begin{Definition}
Let $\sT$ be a triangulated category with set indexed coproducts and
let $B$ be an object of $\sT$.

Then $\langle B \rangle$ denotes the triangulated subcategory of $\sT$
consisting of objects built from $B$ using distinguished triangles,
retracts, and set indexed coproducts (cf.\ \cite[def.\
3.2.9]{Neemanbook}).
\hfill $\Box$
\end{Definition}

\begin{Definition}
Let $\sT$ be a triangulated category with set indexed coproducts.  An
object $B$ of $\sT$ will be called self-compact if the restricted
functor 
\[
  \Hom_{\ssT}(B,-)\big|_{\langle B \rangle}
\]
respects set indexed coproducts.
\hfill $\Box$
\end{Definition}

\begin{Remark}
\label{rmk:Z2a}
Recall that an object $C$ of $\sT$ is called compact if the functor
\[
  \Hom_{\ssT}(C,-)
\]
respects set indexed coproducts.  

A compact object is self-compact, but there are self-compact objects
which are not compact.  For instance, if $\BZ_2$, the integers with
$2$ inverted, are viewed as a complex of $\BZ$-modules, then they are
self-compact but not compact in $\D(\BZ)$, the derived category of the
integers; see example
\ref{exa:Z2b}.
\hfill $\Box$
\end{Remark}

\begin{Remark}
\label{rmk:RHom}
Let $R$ be a DGA with derived category $\D(R)$.  It is not hard to
see that if $B$ is self-compact and $C$ compact in $\D(R)$, then the
functors 
\[
  \RHom_R(B,-)\big|_{\langle B \rangle}
\]
and
\[
  \RHom_R(C,-)
\]
respect set indexed coproducts.
\hfill $\Box$
\end{Remark}

For the following results, recall that if $\cE$ is a DGA then
$\cE^{\opp}$ denotes the opposite DGA with product $\cdot$ defined in
terms of the product of $\cE$ by $e \cdot f = (-1)^{|e||f|}fe$.  Left
DG modules over $\cE^{\opp}$ can be identified canonically with right
DG modules over $\cE$, so $\D(\cE^{\opp})$, the derived category of
left DG modules over $\cE^{\opp}$, can be identified with the derived
category of right DG modules over $\cE$.  Note that subscripts
indicate left and right DG module structures.

\begin{Lemma}
\label{lem:Thomason}
If $\cE$ is a DGA then $\D(\cE^{\opp}) = \langle \cE_{\cE} \rangle$.
\end{Lemma}

\begin{proof}
This is a consequence of Neeman-Thomason localization, \cite[thm.\
2.1.2]{Neeman}.
\end{proof}

\begin{Theorem}
\label{thm:B}
Let $R$ be a DGA with a K-projective left DG module $B$ which is
self-compact in $\D(R)$, and let $\cE$ be the endomorphism DGA of
${}_{R}B$.  Then $B$ acquires the structure ${}_{R,\cE}B$, and there
is an adjoint pair of functors
\[
\xymatrix
{
  \D(\cE^{\opp})
    \ar[rrr]<1ex>^-{i_{*}(-) = - \LTensor_{\cE} B} & & &
  \D(R)
    \ar[lll]<1ex>^-{i^!(-) = \RHom_R(B,-)}
}
\]
where $i_{*}$ is a full embedding with essential image 
\[
  \EssIm i_* = \langle {}_{R}B \rangle.
\]
\end{Theorem}

\begin{proof}
It is clear that $B$ acquires the structure ${}_{R,\cE}B$ since $\cE$
is the endomorphism DGA of ${}_{R}B$, and hence by definition acts on
$B$ in a way compatible with the action of $R$.

Lemma \ref{lem:Thomason} says $\sD(\cE^{\opp}) = \langle \cE_{\cE}
\rangle$, that is, each object in $\sD(\cE^{\opp})$ is built from
$\cE_{\cE}$ using the operations of distinguished triangles, retracts,
and set indexed coproducts.  The functor 
\[
  i_*(-) = - \LTensor_{\cE} {}_{R,\cE}B
\]
respects these operations, so each object in the essential image of
$i_*$ is built from $i_*(\cE_{\cE}) = \cE_{\cE} \LTensor_{\cE}
{}_{R,\cE}B \cong {}_{R}B$ using distinguished triangles, retracts,
and set indexed coproducts; that is,
\begin{equation}
\label{equ:q}
  \EssIm i_* \subseteq \langle {}_{R}B \rangle.
\end{equation}
Since $i_*$ respects set indexed coproducts, it follows that $i_*$
sends set indexed coproducts in $\sD(\cE^{\opp})$ to set indexed
coproducts in $\langle {}_{R}B \rangle$.

Moreover, ${}_{R}B$ is self-compact so the restriction of 
\[
  i^!(-) = \RHom_R({}_{R,\cE}B,-) 
\]
to $\langle {}_{R}B \rangle$ respects set indexed coproducts by remark
\ref{rmk:RHom}.  Together, this shows that the functor 
\[
  i^! i_*(-)
\]
respects set indexed coproducts.

Note that the unit morphism
\[
  \cE_{\cE} \longrightarrow  i^! i_*(\cE_{\cE})
\]
is just the canonical morphism
\[
  \cE_{\cE}
  \longrightarrow
  \RHom_R({}_{R,\cE}B,\cE_{\cE} \LTensor_{\cE} {}_{R,\cE}B)
  \cong \RHom_R({}_{R,\cE}B,{}_{R}B)
\]
which is an isomorphism since $\cE$ is the endomorphism DGA of the
K-projective DG module ${}_{R}B$.  Since $i^! i_*$ respects set
indexed coproducts, it follows that the unit morphism
\[
  Y \longrightarrow i^! i_* Y
\]
is an isomorphism for each object $Y$ which can be built from
$\cE_{\cE}$, that is, for each object in $\langle \cE_{\cE} \rangle =
\D(\cE^{\opp})$. 

By adjoint functor theory this implies that $i_*$ is a full embedding
of $\D(\cE^{\opp})$ into $\D(R)$.  

To conclude the proof, I must show $\EssIm i_* = \langle {}_{R}B
\rangle$.  The inclusion $\subseteq$ was proved in equation
\eqref{equ:q}, so I must show
\[
  \EssIm i_* \supseteq \langle {}_{R}B \rangle.
\]
For this, note that $i_*(\cE_{\cE}) \cong {}_{R}B$ is in the essential
image of $i_*$.  Since $i_*$ is a full embedding respecting set
indexed coproducts, it follows that each object built from ${}_{R}B$
is in the essential image of $i_*$, as desired.
\end{proof}

\begin{Example}
\label{exa:Z2b}
The purpose of this example is to prove that $\BZ_2$, the integers
with $2$ inverted, is self-compact but not compact in $\D(\BZ)$, the
derived category of the integers.

There is an adjoint pair of functors
\[
\xymatrix
{
  \D(\BZ_2)
    \ar[rr]<-1ex>_-{i_*} & &
  \D(\BZ)
    \ar[ll]<-1ex>_-{i^*}
}
\]
where $i^*(-) = \BZ_2 \LTensor_{\BZ} -$ while $i_*$ is the forgetful
functor which takes a complex of $\BZ_2$-modules and views it as a
complex of $\BZ$-modules.

By adjoint functor theory, the unit morphism
\[
  i_* \BZ_2 \longrightarrow i_*i^*(i_* \BZ_2)
\]
is an isomorphism.  The functors $i_*$ and $i^*$ clearly respect set
indexed coproducts, so $i_* i^*$ respects set indexed coproducts, and
it follows that the unit morphism
\[
  X \longrightarrow i_*i^*X
\]
is an isomorphism for each $X$ in $\langle i_* \BZ_2 \rangle$.

This permits the computation
\begin{align*}
  \Hom_{\Dsmall(\BZ)}(i_* \BZ_2,-)\big|_{\langle i_* \BZ_2 \rangle}
  & \simeq \Hom_{\Dsmall(\BZ)}
           (i_* \BZ_2,i_*i^*(-))\big|_{\langle i_* \BZ_2 \rangle} \\
  & \simeq
           \Hom_{\Dsmall(\BZ_2)}
           (i^*i_* \BZ_2,i^*(-))\big|_{\langle i_* \BZ_2 \rangle} \\
  & \stackrel{\rm (a)}{\simeq}
           \Hom_{\Dsmall(\BZ_2)}
           (\BZ_2,i^*(-))\big|_{\langle i_* \BZ_2 \rangle},
\end{align*}
where (a) is because $i^*i_*\BZ_2 = \BZ_2 \LTensor_{\BZ} i_* \BZ_2 =
\BZ_2 \LTensor_{\BZ} \BZ_2 \cong \BZ_2$.  But $i^*$ respects set
indexed coproducts, so the same holds for the right-hand side of the
computation and in consequence for the left-hand side,
\[
  \Hom_{\Dsmall(\BZ)}(i_* \BZ_2,-)\big|_{\langle i_* \BZ_2 \rangle}.
\]
This shows that $i_* \BZ_2$ is a self-compact object of $\D(\BZ)$, and
$i_* \BZ_2$ is just $\BZ_2$ viewed as a complex of $\BZ$-modules, so
$\BZ_2$ is self-compact in $\D(\BZ)$.

On the other hand, $\BZ_2$ is not finitely generated over $\BZ$ so
cannot be compact in $\D(\BZ)$ since the compact objects in $\D(\BZ)$
have finitely generated cohomology as follows from \cite[thm.\
5.3]{KellerDG}.
\hfill $\Box$
\end{Example}

\section{Dwyer and Greenlees's Morita theory}
\label{sec:DwyerGreenlees}

\begin{Setup}
\label{set:DwyerGreenleesSetup}
The following is taken from \cite{DwyerGreenlees}, up to the
trivial change of $R$ being a DGA and not a ring.

Let $R$ be a DGA with a K-projective left DG module $C$ which is
compact in $\D(R)$, and let $\cF$ be the endomorphism DGA of
${}_{R}C$.  Then $C$ acquires the structure ${}_{R,\cF}C$.

Observe that
\begin{align*}
  \RHom_R({}_{R,\cF}C,-) 
  & \simeq
    \RHom_R({}_{R,\cF}C,{}_{R}R_{R} \LTensor_R -) \\
  & \stackrel{\rm (a)}{\simeq}
    \RHom_R({}_{R,\cF}C,{}_{R}R_{R}) \LTensor_R -
\end{align*}
where (a) is because ${}_{R}C$ is compact.  Setting
\[
  C^*_{R,\cF} = \RHom_R({}_{R,\cF}C,{}_{R}R_{R})
\]
hence gives
\[
  \RHom_R({}_{R,\cF}C,-)
  \simeq C^*_{R,\cF} \LTensor_R -.
\]

There are therefore functors
\begin{equation}
\label{equ:e}
\xymatrix
{
&&\D(R)
    \ar[rr]^{j^*} & &
  \D(\cF^{\opp})
    \ar@/_1.5pc/[ll]_{j_!} \ar@/^1.5pc/[ll]^{j_*} & &
}
\end{equation}
given by
\begin{align*}
  j_!(-) & = - \LTensor_{\cF} {}_{R,\cF}C, \\
  j^*(-) & = \RHom_R({}_{R,\cF}C,-) \simeq C^*_{R,\cF} \LTensor_R -, \\
  j_*(-) & = \RHom_{\cF^{\opp}}(C^*_{R,\cF},-),
\end{align*}
where $(j_!,j^*)$, and $(j^*,j_*)$ are adjoint pairs.  
\hfill $\Box$
\end{Setup}

The following result was established in \cite[sec.\
2]{DwyerGreenlees}, up to the change of $R$ being a DGA and not a
ring.

\begin{Proposition}
\label{pro:DwyerGreenlees}
In the situation of setup \ref{set:DwyerGreenleesSetup}, the functors
$j_!$ and $j_*$ are full embeddings.
\end{Proposition}

\section{Recollement}
\label{sec:recollement}

Let me first recall the definition of recollement from
\cite[sec.\ 1.4]{BBD}. 

\begin{Definition}
\label{dfn:recollement}
A recollement of triangulated categories is a diagram of triangulated
categories and triangulated functors
\[
\xymatrix
{
&&\sT^{\prime} 
    \ar[rr]^{i_*} & & 
  \sT 
    \ar[rr]^{j^*}
    \ar@/_1.5pc/[ll]_{i^*} \ar@/^1.5pc/[ll]^{i^!} & &
  \sT^{\prime\prime}
    \ar@/_1.5pc/[ll]_{j_!} \ar@/^1.5pc/[ll]^{j_*} & &
}
\]
satisfying
\begin{enumerate}

  \item  $(i^*,i_*)$, $(i_*,i^!)$, $(j_!,j^*)$, and $(j^*,j_*)$ are
         adjoint pairs.

\smallskip

  \item  $j^* i_* = 0$.

\smallskip

  \item  $i_*$, $j_!$, and $j_*$ are full embeddings.

\smallskip

  \item  Each object $X$ in $\sT$ determines distinguished triangles
\smallskip
\begin{enumerate}

  \item  $i_*i^!X \longrightarrow X \longrightarrow j_*j^*X
         \longrightarrow \;\;$ and

  \item  $j_!j^*X \longrightarrow X \longrightarrow i_*i^*X
         \longrightarrow$ 

\end{enumerate}
\smallskip
where the arrows to and from $X$ are counit and unit morphisms.
\end{enumerate}
\hfill $\Box$
\end{Definition}

\begin{Remark}
\label{rmk:recollement}
The following are easy formal consequences of definition
\ref{dfn:recollement}.
\begin{enumerate}

  \item  $i^* j_! = 0$ and $i^! j_* = 0$.
%
%

\smallskip

  \item  The restriction of $i_*i^*$ to the essential image of $i_*$
         is naturally equivalent to the identity functor.

\end{enumerate}
\hfill $\Box$
\end{Remark}

For the following results, note that if $X$ is a full subcategory of a
triangulated category $\sT$, then there are full subcategories
\[
  X^{\perp} 
  = \{\, Y \in \sT \,|\, 
         \Hom_{\ssT}(\Sigma^{\ell} X,Y) = 0 \mbox{ for each } \ell \,\}
\]
and
\[
  {}^{\perp}X
  = \{\, Y \in \sT \,|\, 
         \Hom_{\ssT}(Y,\Sigma^{\ell} X) = 0 \mbox{ for each } \ell \,\}.
\]

It turns out that Dwyer and Greenlees's Morita theory can be improved
in terms of recollement.  Specifically, equation \eqref{equ:e} from
setup \ref{set:DwyerGreenleesSetup} is the right-hand part of a
recollement as follows.

\begin{Proposition}
\label{pro:DwyerGreenleesimproved}
In the situation of setup \ref{set:DwyerGreenleesSetup}, there is a
recollement
\[
\xymatrix
{
&&({}_{R}C)^{\perp}
    \ar[rr]^{i_*} & & 
  \D(R)
    \ar[rr]^{j^*}
    \ar@/_1.5pc/[ll]_{i^*} \ar@/^1.5pc/[ll]^{i^!} & &
  \D(\cF^{\opp})
    \ar@/_1.5pc/[ll]_{j_!} \ar@/^1.5pc/[ll]^{j_*} & &
}
\]
where $i_*$ is the inclusion of the full subcategory
$({}_{R}C)^{\perp}$ and $i^*$ and $i^!$ are its left- and
right-adjoint functors, while the functors $j_!$, $j^*$, and $j_*$ are
given as in setup \ref{set:DwyerGreenleesSetup},
\begin{align*}
  j_!(-) & = - \LTensor_{\cF} {}_{R,\cF}C, \\
  j^*(-) & = \RHom_R({}_{R,\cF}C,-), \\
  j_*(-) & = \RHom_{\cF^{\opp}}(C^*_{R,\cF},-).
\end{align*}
\end{Proposition}

\begin{proof}
The functors $j_!$ and $j_*$ are left- and right-adjoint to $j^*$, and
by proposition \ref{pro:DwyerGreenlees} both $j_!$ and $j_*$ are full
embeddings.

This situation is exactly the one considered in \cite[prop.\
2.7]{Miyachi} which now gives existence of a recollement where the
left-hand category is the kernel of $j^*$ and $i_*$ is the inclusion.
But the kernel of $j^*(-) = \RHom_R({}_{R,\cF}C,-)$ is clearly
$({}_{R}C)^{\perp}$, so the present proposition follows.
\end{proof}

The preceding material allows me to prove the following main result.

\begin{Theorem}
\label{thm:main}
Let $R$ be a DGA with left DG modules $B$ and $C$.  Then the following
are equivalent.
\begin{enumerate}

  \item  There is a recollement
\[
\xymatrix
{
&&\D(S)
    \ar[rr]^{i_*} & & 
  \D(R)
    \ar[rr]^{j^*}
    \ar@/_1.5pc/[ll]_{i^*} \ar@/^1.5pc/[ll]^{i^!} & &
  \D(T)
    \ar@/_1.5pc/[ll]_{j_!} \ar@/^1.5pc/[ll]^{j_*} & &
}
\]
where $S$ and $T$ are DGAs, for which
\[
  i_*({}_{S}S) \cong B, \;\; j_!({}_{T}T) \cong C.
\]

\smallskip

  \item  In the derived category $\D(R)$, the DG module $B$ is
         self-compact, $C$ is compact, $B^{\perp} \cap C^{\perp} = 0$,
         and $B \in C^{\perp}$.

\end{enumerate}
\end{Theorem}

\begin{proof}
(i) $\Rightarrow$ (ii) \hspace{0.5ex}
The functor $i_*$ is triangulated and a full embedding by definition
\ref{dfn:recollement}(iii); hence the essential image of $i_*$ is a
triangulated subcategory of $\D(R)$.  And $i_*$ is a left-adjoint so
respects set indexed coproducts; hence the essential image of $i_*$
is closed under set indexed coproducts in $\D(R)$.

Since $B \cong i_*({}_{S}S)$ is in the essential image of $i_*$, it
follows that 
\begin{equation}
\label{equ:k}
  \langle B \rangle \subseteq \EssIm i_*,
\end{equation}
and remark \ref{rmk:recollement}(ii) then implies that the restriction
of $i_*i^*$ to $\langle B \rangle$ is naturally equivalent to the
identity functor.  This permits the computation
\begin{align*}
  \Hom_{\Dsmall(R)}(B,-) \big|_{\langle B \rangle}
  & \simeq \Hom_{\Dsmall(R)}(i_*({}_{S}S),i_*i^*(-))
              \big|_{\langle B \rangle} \\
  & \stackrel{\rm (a)}{\simeq} \Hom_{\Dsmall(S)}({}_{S}S,i^*(-))
              \big|_{\langle B \rangle},
\end{align*}
where (a) is because $i_*$ is a full embedding.  But
$i^*$ is a left-adjoint so respects set indexed coproducts, so the
same holds for the right-hand side of the computation and in
consequence for the left-hand side,
\[
  \Hom_{\Dsmall(R)}(B,-) \big|_{\langle B \rangle}.
\]
This shows that $B$ is self-compact.

Similarly, there is the computation
\[
  \Hom_{\Dsmall(R)}(C,-)
    \simeq \Hom_{\Dsmall(R)}(j_!({}_{T}T),-) 
    \simeq \Hom_{\Dsmall(T)}({}_{T}T,j^*(-)),
\]
and $j^*$ is a left-adjoint so respects set indexed coproducts, so
the same holds for the right-hand side of the computation and in
consequence for the left-hand side,
\[
  \Hom_{\Dsmall(R)}(C,-).
\]
This shows that $C$ is compact.

Let $X$ be in $B^{\perp} \cap C^{\perp}$.  Then
\begin{align*}
  0 & = \Hom_{\Dsmall(R)}(\Sigma^{\ell} B,X) \\
    & \cong \Hom_{\Dsmall(R)}(\Sigma^{\ell} i_*({}_{S}S),X) \\
    & \cong \Hom_{\Dsmall(S)}(\Sigma^{\ell}({}_{S}S),i^! X)
\end{align*}
and
\begin{align*}
  0 & = \Hom_{\Dsmall(R)}(\Sigma^{\ell} C,X) \\
    & \cong \Hom_{\Dsmall(R)}(\Sigma^{\ell} j_!({}_{T}T),X) \\
    & \cong \Hom_{\Dsmall(T)}(\Sigma^{\ell}({}_{T}T),j^* X)
\end{align*}
for each $\ell$, proving $i^! X = 0 = j^* X$.  But then the
distinguished triangle in definition \ref{dfn:recollement}(iv)(a)
shows $X = 0$, and $B^{\perp} \cap C^{\perp} = 0$ follows.

Finally,
\begin{align*}
  \Hom_{\Dsmall(R)}(\Sigma^{\ell}C,B)
  & \cong
    \Hom_{\Dsmall(R)}(\Sigma^{\ell}j_!({}_{T}T),i_*({}_{S}S)) \\
  & \cong
    \Hom_{\Dsmall(T)}(\Sigma^{\ell}({}_{T}T),j^*i_*({}_{S}S)),
\end{align*}
and this is $0$ for each $\ell$ because $j^*i_* = 0$ by definition
\ref{dfn:recollement}(ii), so $B$ is in $C^{\perp}$.

\medskip
\noindent
(ii) $\Rightarrow$ (i) \hspace{0.5ex} 
It is enough to construct a recollement
\begin{equation}
\label{equ:g}
\hspace{-69ex}
\lefteqn{
\xymatrix
{
&&\D(\cE^{\opp})
    \ar[rr]^{i_*} & & 
  \D(R)
    \ar[rr]^{j^*}
    \ar@/_1.5pc/[ll]_{i^*} \ar@/^1.5pc/[ll]^{i^!} & &
  \D(\cF^{\opp})
    \ar@/_1.5pc/[ll]_{j_!} \ar@/^1.5pc/[ll]^{j_*} & &
}
}
\end{equation}
for which
\begin{equation}
\label{equ:h}
  i_*(\cE_{\cE}) \cong B, \;\; j_!(\cF_{\cF}) \cong C,
\end{equation}
because the recollement in part (i) of the theorem can be obtained
from this by setting $S = \cE^{\opp}$ and $T = \cF^{\opp}$.

I can clearly replace $B$ and $C$ with K-projective resolutions.  Let
$\cE$ and $\cF$ be the endomorphism DGAs of ${}_{R}B$ and ${}_{R}C$
so I have the full embedding of theorem \ref{thm:B} because ${}_{R}B$
is self-compact, and the recollement of proposition
\ref{pro:DwyerGreenleesimproved} because ${}_{R}C$ is compact.

The recollement of proposition \ref{pro:DwyerGreenleesimproved} goes
some way towards giving \eqref{equ:g}, except that the left-hand
category is $({}_{R}C)^{\perp}$ and not $\D(\cE^{\opp})$.  But if I
could prove
\begin{equation}
\label{equ:j}
  ({}_{R}C)^{\perp} = \langle {}_{R}B \rangle,
\end{equation}
then I could replace $({}_{R}C)^{\perp}$ by $\langle {}_{R}B
\rangle$ which could again be replaced by $\D(\cE^{\opp})$ using
the full embedding of theorem \ref{thm:B}, and this would give
\eqref{equ:g}.  In this case, \eqref{equ:h} would be clear because
theorem \ref{thm:B} would imply
\[
  i_*(\cE_{\cE}) = \cE_{\cE} \LTensor_{\cE} {}_{R,\cE}B
  \cong {}_{R}B
\]
while proposition \ref{pro:DwyerGreenleesimproved} would imply
\[
  j_!(\cF_{\cF}) = \cF_{\cF} \LTensor_{\cF} {}_{R,\cF}C
  \cong {}_{R}C.
\]

To show \eqref{equ:j}, note that $\supseteq$ is clear since ${}_{R}B$
is in $({}_{R}C)^{\perp}$ by assumption while ${}_{R}C$ is compact.  To prove
$\subseteq$, let $X$ be in $({}_{R}C)^{\perp}$.  The adjunction in
theorem \ref{thm:B} gives a counit morphism $i_*i^!X
\stackrel{\epsilon}{\longrightarrow} X$ (where $i_*$ and $i^!$ are now
used in the sense of theorem \ref{thm:B}) which can be extended to a
distinguished triangle 
\[
  i_*i^!X 
  \stackrel{\epsilon}{\longrightarrow} X
  \longrightarrow Y \longrightarrow.
\]
By adjoint functor theory, $i^!(\epsilon)$ is an isomorphism, so $i^!Y
= 0$, that is, $\RHom_R({}_{R,\cE}B,Y) = 0$, so $Y$ is in
$({}_{R}B)^{\perp}$. 

Moreover, $i_*i^!X$ is in the essential image of $i_*$ which equals
$\langle {}_{R}B \rangle$ by theorem \ref{thm:B}, so since $\langle
{}_{R}B \rangle \subseteq ({}_{R}C)^{\perp}$ it follows that $i_*i^!X$
is in $({}_{R}C)^{\perp}$.  But $X$ is in $({}_{R}C)^{\perp}$ by
assumption, and it follows that also $Y$ is in $({}_{R}C)^{\perp}$.

So $Y$ is in $({}_{R}B)^{\perp} \cap ({}_{R}C)^{\perp}$ which is $0$
by assumption, so $Y = 0$, so the distinguished triangle shows $X
\cong i_*i^!X$ and this is in the essential image of $i_*$ which is
equal to $\langle {}_{R}B \rangle$.
\end{proof}

\begin{Remark}
\label{rmk:main}
Note that the proof of theorem \ref{thm:main}, (ii) $\Rightarrow$ (i),
gives a recipe for constructing $S$ and $T$ when $R$, $B$, and $C$ are
known:

Replace $B$ and $C$ with K-projective resolutions, set $\cE$ and $\cF$
equal to the endomorphism DGAs of ${}_{R}B$ and ${}_{R}C$, and set $S
= \cE^{\opp}$ and $T = \cF^{\opp}$.

Similarly, there is a recipe for constructing the functors $i_*$,
$i^!$, $j_!$, $j^*$, and $j_*$:

After replacing $B$ and $C$ with K-projective resolutions, $B$ and $C$
acquire the structures ${}_{R,\cE}B$ and ${}_{R,\cF}C$, that is,
${}_{R}B_{S}$ and ${}_{R}C_{T}$, and the functors are then given by
\begin{eqnarray*}
                                    & & j_!(-) = {}_{R}C_{T} \LTensor_T - \\
  i_*(-) = {}_{R}B_{S} \LTensor_S - & & j^*(-) = \RHom_R({}_{R}C_{T},-) \\
  i^!(-) = \RHom_R({}_{R}B_{S},-)   & & j_*(-) = \RHom_{T}({}_{T}C^*_{R},-).
\end{eqnarray*}
\hfill $\Box$
\end{Remark}

\begin{Example}
Let 
\[
  R = \BZ, \; B = \BZ_2, \; C = \BZ/(2)
\]
where $\BZ_2$ is $\BZ$ with $2$ inverted.  The purpose of this example
is to show that these data satisfy the conditions of theorem
\ref{thm:main}(ii).  Since example \ref{exa:Z2b} proves that $\BZ_2$
is self-compact but not compact in $\D(\BZ)$, this shows that
theorem \ref{thm:main} really needs the notion of
self-com\-pact\-ness.

To check the conditions of theorem \ref{thm:main}(ii) apart from
self-com\-pact\-ness of $B = \BZ_2$ which is already known, note that
there is a distinguished triangle
\[
  \BZ \stackrel{2}{\longrightarrow} \BZ \longrightarrow \BZ/(2) \longrightarrow
\]
in $\D(\BZ)$.  So $\BZ/(2)$ is finitely built from $\BZ$ and hence
$\BZ/(2)$ is compact in $\D(\BZ)$; that is, $C$ is compact in $\D(R)$.

Also, the distinguished triangle gives a distinguished triangle
\[
  \RHom_{\BZ}(\BZ/(2),X) \longrightarrow
  \RHom_{\BZ}(\BZ,X) \stackrel{2}{\longrightarrow}
  \RHom_{\BZ}(\BZ,X) \longrightarrow
\]
for each $X$ in $\D(\BZ)$, that is,
\[
  \RHom_{\BZ}(\BZ/(2),X) \longrightarrow
  X \stackrel{2}{\longrightarrow}
  X \longrightarrow,
\]
and the long exact sequence of this implies that
$\RHom_{\BZ}(\BZ/(2),X) = 0$ if and only if $2$ acts invertibly on
each cohomology module of $X$.

So for example, $\RHom_{\BZ}(\BZ/(2),\BZ_2) = 0$, and this implies
that $\BZ_2$ is in $(\BZ/(2))^{\perp}$, that is, $B$ is in
$C^{\perp}$.

Finally, let $X$ be in $(\BZ_2)^{\perp} \cap (\BZ/(2))^{\perp}$.  It
is well known that since $\BZ$ has global dimension one, $X$ is
isomorphic in $\D(\BZ)$ to the complex $\widetilde{X}$ having
$\H^i(X)$ in cohomological degree $i$ and having zero differential.
It was shown above that since $X$ is in $(\BZ/(2))^{\perp}$, the
integer $2$ acts invertibly on each $\H^i(X)$.  That is, each
$\H^i(X)$ is in fact a $\BZ_2$-module, so $\widetilde{X}$ can be
viewed as a complex of $\BZ_2$-modules which I will denote $Y$.

Now the forgetful functor $i_*$ from example \ref{exa:Z2b} satisfies
\begin{equation}
\label{equ:l}
  i_* Y = \widetilde{X} \cong X, 
\end{equation}
and using also the left-adjoint functor $i^*$ from example
\ref{exa:Z2b}, I can compute for each $\ell$,
\begin{align*}
  \Hom_{\Dsmall(\BZ_2)}(\Sigma^{\ell} \BZ_2,Y)
  & \stackrel{\rm (a)}{\cong} 
    \Hom_{\Dsmall(\BZ_2)}(\Sigma^{\ell} i^* \BZ_2,Y) \\
  & \cong \Hom_{\Dsmall(\BZ)}(\Sigma^{\ell}\BZ_2,i_* Y) \\
  & \stackrel{\rm (b)}{\cong}
    \Hom_{\Dsmall(\BZ)}(\Sigma^{\ell}\BZ_2,X) \\
  & \stackrel{\rm (c)}{\cong} 0,
\end{align*}
where (a) is because $\BZ_2 \cong \BZ_2 \LTensor_{\BZ} \BZ_2
\cong i^*\BZ_2$ while (b) is by equation \eqref{equ:l} and (c) is
because $X$ is in $(\BZ_2)^{\perp}$.

This proves $Y = 0$ whence $X = 0$ by equation \eqref{equ:l}, and
altogether, I have shown $(\BZ_2)^{\perp} \cap (\BZ/(2))^{\perp} = 0$,
that is, $B^{\perp} \cap C^{\perp} = 0$.
\hfill $\Box$
\end{Example}

The proof of the following theorem could have used the concrete
construction of a recollement at the end of the proof of theorem
\ref{thm:main} as a crutch.  However, a direct proof is just
as easy.

\begin{Theorem}
\label{thm:main2}
In the situation of theorem \ref{thm:main}, the following hold.
\begin{enumerate}

  \item  $\EssIm i_* = \Ker j^* = \langle B \rangle = C^{\perp} =
         {}^{\perp}(B^{\perp})$.

\smallskip

  \item  $\EssIm j_! = \Ker i^* = \langle C \rangle =
         {}^{\perp}(C^{\perp})$. 

\smallskip

  \item  $\EssIm j_* = \Ker i^! = B^{\perp} = (C^{\perp})^{\perp}$.

\end{enumerate}
\end{Theorem}

\begin{proof}
To prove $\EssIm i_* = \Ker j^*$ in (i), note that the inclusion
$\subseteq$ holds because $j^*i_* = 0$ by definition
\ref{dfn:recollement}(ii), and that the inclusion $\supseteq$ follows
from the distinguished triangle in definition
\ref{dfn:recollement}(iv)(b).

The equalities $\EssIm j_! = \Ker i^*$ and $\EssIm j_* = \Ker i^!$ in
(ii) and (iii) are proved by similar arguments.

To prove $\EssIm i_* = \langle B \rangle$ in (i), note that the
inclusion $\supseteq$ is already known from equation
\eqref{equ:k} in the proof of theorem \ref{thm:main}.  To see the
inclusion $\subseteq$, it is enough to see $i_*(\D(S)) \subseteq
\langle B \rangle$.  But this is easy,
\[
  i_*(\D(S)) 
  \stackrel{\rm (a)}{=} i_*(\langle {}_{S}S \rangle)
  \stackrel{\rm (b)}{\subseteq} \langle i_*({}_{S}S) \rangle 
  = \langle B \rangle,
\]
where (a) is by lemma \ref{lem:Thomason} and (b) is because $i_*$ is a
left-adjoint so respects set indexed coproducts.

The equality $\EssIm j_! = \langle C \rangle$ in (ii) is proved by a
similar argument.

The equality $\langle B \rangle = C^{\perp}$ in (i) is already known
from equation \eqref{equ:j} in the proof of theorem \ref{thm:main}.

To show $\Ker i^* = {}^{\perp}(C^{\perp})$ in (ii), note that 
$X$ is in $\Ker i^*$ if and only if
\begin{equation}
\label{equ:p}
  0 = \Hom_{\Dsmall(S)}(i^* X,Y)
    \cong \Hom_{\Dsmall(R)}(X,i_* Y)
\end{equation}
for each $Y$ in $\D(S)$.  But I have already proved $\EssIm i_* =
C^{\perp}$, so up to isomorphism, the objects of the form $i_* Y$ are
exactly the objects in $C^{\perp}$, so equation
\eqref{equ:p} is equivalent to $X$ being in ${}^{\perp}(C^{\perp})$.

Similar arguments show first the equality $\Ker i^! = B^{\perp}$ in
(iii) and then $\Ker j^* = {}^{\perp}(B^{\perp})$ in (i).

Finally, $B^{\perp} = (C^{\perp})^{\perp}$ in (iii) follows from
$\langle B \rangle = C^{\perp}$ which I have already proved. 
\end{proof}

\medskip
\noindent
{\bf Acknowledgement. }  The diagrams were typeset with XY-pic.

\end{document}